\documentclass[final,3p, times, 12pt]{elsarticle}

\usepackage{amssymb,amsmath,color}
\usepackage[mathscr]{eucal}
 \usepackage{lineno}

 \journal{ }
\renewcommand{\journal}{  TBD }
\newtheorem{thm}{Theorem}

\newtheorem{example}{Example}
\newtheorem{defin}{Definition}
\newproof{pf}{Proof}

\begin{document}

\begin{frontmatter}



\title{A construction of pooling designs with surprisingly high degree of error correction}

\author[JG]{Jun Guo}\ead{guojun$_-$lf@163.com}
\author[KW]{Kaishun Wang\corref{cor}}
\ead{wangks@bnu.edu.cn}
\cortext[cor]{Corresponding author}

\address[JG]{Math. and Inf. College, Langfang Teachers'
College, Langfang  065000,  China }
\address[KW]{Sch. Math. Sci. \& Lab. Math. Com. Sys.,
Beijing Normal University, Beijing  100875, China}

\begin{abstract}
It is well-known that many famous pooling designs are constructed
from mathematical structures by the  ``containment matrix" method.
  In this paper, we propose another method   and
obtain a family of pooling designs with surprisingly high degree of
error correction based on a finite set. Given  the numbers of items
and pools, the error-tolerant property of our designs is much better
than that of Macula's designs when the size of the set is large
enough.

\end{abstract}

\end{frontmatter}

\section*{}

Pooling design is a mathematical tool to reduce the number of tests
in DNA library screening \cite{CD,CD2,DH}. A pooling design is
usually represented by a binary matrix with columns  indexed with
items and rows
 indexed with pools. A cell $(i,j)$ contains a 1-entry if and
only if the $i$th pool contains the $j$th item. Biological
experiments are notorious for producing erronous outcomes.
Therefore, it would be wise for pooling designs to allow some
outcomes to be affected by errors.  A binary matrix $M$ is called
$s^e$-{\em disjunct} if given any $s+1$ columns of $M$ with one
designated, there are $e+1$ rows with a 1 in the designated column
and 0 in each of the other $s$ columns. A $s^0$-disjunct matrix   is
also called   $s$-{\em disjunct}. An $s^e$-disjunct matrix is called
{\it fully $s^e$-disjunct} if it is not $s_1^{e_1}$-disjunct
whenever $s_1>s$ or $e_1>e$. An $s^e$-disjunct matrix is $\lfloor
e/2\rfloor$-error-correcting (see \cite{DF}).

For positive integers $k\leq n$, let $[n]=\{1,2,\ldots,n\}$ and
$\left([n]\atop k\right)$ be the set of all $k$-subsets of
$[n]$.

Macula \cite{Macula,Macula2} proposed a novel way of constructing
disjunct matrices by the containment relation of subsets in a finite
set.

\begin{defin}{\rm(\cite{Macula})}\,  For positive integers $1\leq d<k< n$,
let $M(d,k,n)$ be the binary
matrix with rows indexed with $\left([n]\atop
d\right)$ and columns indexed with  $\left([n]\atop k\right)$ such that $M(A,B)=1$ if
and only if $A\subseteq B$.
\end{defin}

D'yachkov et al. \cite{DF2} discussed the error-correcting property
of $M(d,k,n)$.

\begin{thm}{\rm (\cite{DF2})}\label{thmdy} For positive integers $1\leq d<k< n$ and $s\leq d$, $M(d,k,n)$ is fully
$s^{e_1}$-disjunct,   where $e_1=\left(k-s\atop d-s\right)-1$.
\end{thm}

Ngo and Du \cite{Ngo} constructed disjunct matrices by the
containment relation of subspaces in a finite vector space.
D'yachkov et al. \cite{DF} discussed the error-tolerant property of
Ngo and Du's construction. Huang and Weng \cite{HW} introduced the
comprehensive concept of pooling spaces, which is a significant
addition to the general theory. Recently,
 many pooling designs have been constructed using the ``containment
matrix" method, see e.g. \cite{bai, hhw, hww}.

Next we shall introduce our construction.

\begin{defin}
Given    integers $1\leq d< k<n $ and $0\leq i\leq
 d$.  Let $M(i;d,k,n)$ be the binary matrix with rows indexed
$\left([n]\atop d\right)$ and columns indexed with
$\left([n]\atop k\right)$ such that $M(A,B)=1$ if and only if
$|A\cap B|=i$.
\end{defin}

Note that $M(i;d,k,n)$ and $M(d,k,n)$ have the same size, and
$M(i;d,k,n)$ is an $\left(n\atop d\right)\times \left(n\atop
k\right)$ matrix with row weight $\left(d\atop
i\right)\left(n-d\atop k-i\right)$ and column weight $\left(k\atop
i\right)\left(n-k\atop d-i\right)$. Since $M(d;d,k,n)=M(d,k,n)$,
  our construction is a generalization of Macula's matrix.

Let $B\in\left([n]\atop k\right)$ and $C=[n]\backslash B$. Then, for any
$D\in\left([n]\atop d\right)$, $|D\cap B|=i$ if and only if $|D\cap
C|=d-i$. Therefore,
 $M(i;d,k,n)=M(d-i;d,n-k,n)$ when $n>k+d-i$. Since $i\leq\lfloor d/2\rfloor$ if
and only if $d-i\geq\lfloor(d+1)/2\rfloor$,  we always assume that
$i\geq\lfloor(d+1)/2\rfloor$ in this case.

\begin{thm}\label{thm2.1}
Let $1\leq s\leq i,\lfloor(d+1)/2\rfloor\leq i\leq d<k$ and
$n-k-s(k+d-2i)\geq d-i$. Then

\begin{itemize}

\item[\rm(i)] $M(i;d,k,n)$ is an $s^{e_2}$-disjunct matrix, where
$e_2=\left(k-s\atop i-s\right)\left(n-k-s(k+d-2i)\atop
d-i\right)-1$;

\item[\rm(ii)] For a given $k$, if $i<d,$ then $\lim\limits_{n\rightarrow\infty}
\frac{e_2+1}{e_1+1}=\infty.$
\end{itemize}
\end{thm}

\begin{pf} (i)
Let $B_{0},B_{1},\ldots,B_{s}\in\left([n]\atop k\right)$ be any
$s+1$ distinct columns of $M(i;d,k,n)$. Then, for each $j\in[s]$,
there exists an $x_j$ such that $x_j\in B_0\backslash B_j.$ Suppose
$X_0=\{x_j\mid 1\leq j\leq s\}.$ Then $X_0\subseteq B_0$, and
$X_0\not\subseteq B_j$ for each $j\in[s]$.  Note that the number of
$i$-subsets of $B_0$ containing $X_0$ is $\left(k-|X_0|\atop
i-|X_0|\right)= \left(k-|X_0|\atop k-i\right)$. Since $
\left(k-|X_0|\atop k-i\right)$ is decreasing for $1\leq |X_0|\leq s$
and gets its minimum at $|X_0|=s$, the number of $i$-subsets of
$B_0$ containing $X_0$ is at least $\left(k-s\atop k-i\right)$.

Let  $A_0$  be an $i$-subset of $B_0$ containing $X_0$. Then
$|A_0\cap B_j|<i$ for each $j\in[s]$.
 Let $D\in\left([n]\atop
d\right)$ satisfying  $|D\cap B_0|=i$. If there exists $j\in [s]$
such that  $|D\cap B_j|=i$, then $|B_0\cap B_j|\geq|D\cap B_0\cap
B_j|\geq 2i-d$. Suppose   $|B_0\cap B_j|\geq 2i-d$ for each
$j\in[s]$. Since $|\bigcup_{0\leq j\leq s}B_j|\leq k+s(k+d-2i)$, the
number of $d$-subsets $D$ of $[n]$ containing $A_0$  satisfying
$|D\cap B_{0}|=i$ and $|D\cap B_{j}|\not=i$ for each $j\in[s]$ is at
least $\left(n-k-s(k+d-2i)\atop d-i\right)$. Then the number of
$d$-subsets $D$ containing $X_0$ in $\left([n]\atop d\right)$
satisfying $|D\cap B_{0}|=i$ and $|D\cap B_{j}|\not=i$ for each
$j\in[s]$ is at least $\left(k-s\atop
i-s\right)\left(n-k-s(k+d-2i)\atop d-i\right)$. Therefore, (i)
holds.

(ii) is straightforward by (i)  and Theorem~\ref{thmdy}.\qed
\end{pf}

 \begin{example}
$M(5,7,50)$ is fully $1^{14},2^{9}$ and $3^5$-disjunct, but
$M(3;5,7,50)$ is $1^{9989},2^{2324}$ and $3^{299}$-disjunct;
$M(4,5,13)$ is fully $1^3$ and $2^2$-disjunct, but $M(3;4,5,13)$ is
$1^{29}$ and $2^5$-disjunct.
\end{example}

\section*{Concluding remarks}

(i) For  given integers $d<k$ the following limit holds: $\lim\limits_{n\rightarrow
\infty} \frac{\left(n\atop d\right)}{\left(n\atop k\right)}=0$.
This shows that the test-to-item of $M(i;d,k,n)$ is small enough
when $n$ is large enough. By Theorem~\ref{thm2.1}, our pooling
design are better than Macula's designs when $n$ is large enough.

(ii) It seems to be interesting to compute $e$ such  that
$M(i;d,k,n)$ is fully $s^e$-disjunct.

(iii)  In \cite{NG}, Nan and the first author discussed the similar
construction of $s^e$-disjunct matrices in a finite vector space,
but the number $e$ is not well expressed. By the method of this
paper, $e$ may be larger. We will study this problem in a separate
paper.

(iv) For positive integers $1\leq d<k< n$, let $I$ be a nonempty
proper subset of $\{0,1,\ldots, d\}$, and let $M(I;d,k,n)$ be the
binary matrix with rows indexed with
$\left([n]\atop d\right)$ and columns indexed with $\left([n]\atop k\right)$ such
that $M(A,B)=1$ if and only if $|A\cap B|\in I$. How
about the error-tolerant property of $M(I;d,k,n)$?

\section*{Acknowledgment}
We would like thank the referees for their valuable suggestions.
This research is partially supported by    NSF of China,
NCET-08-0052, Langfang Teachers' College (LSZB201005), and   the
Fundamental Research Funds for the Central Universities of China.


\begin{thebibliography}{00}\frenchspacing

\bibitem{bai}

Y. Bai, T. Huang and K. Wang, Error-correcting pooling designs
associated with some distance-regular graphs,   Discrete Appl. Math.
157 (2009) 3038--3045.

\bibitem{CD}
Y. Cheng and D. Du, Efficient constructions of disjunct matrices
with applications to DNA library screening, J. Comput. Biol. 14 (2007) 1208--1216.

\bibitem{CD2}
Y. Cheng and D. Du, New constructions of one- and two-stage pooling sesigns,
J. Comput. Biol. 15 (2008) 195--205.

\bibitem{DH}
D. Du and F.K. Hwang,
 Pooling designs and nonadaptive group testing, Important tools for DNA
 sequencing, Series on Applied Mathematics, 18, World Scientific Publishing Co.
Pte. Ltd., Hackensack, NJ, 2006.

\bibitem{DF}
A. G. D'yachkov, F. K. Hwang, A. J. Macula, P. A. Vilenkin and C. Weng,
A construction of pooling designs with some happy surprises, J. Comput. Biol. 12 (2005) 1127--1134.

\bibitem{DF2}
A. G. D'yachkov, A. J. Macula and P. A. Vilenkin,
Nonadaptive and trivial two-stage group testing with
error-correcting $d^e$-disjunct inclusion matrices, In: Entropy, Search, Complexity,
Bolyai society mathematical studied, vol. 16, Spring, Berlin, pp 71--83, 2007.



\bibitem{hhw} H. Huang, Y. Huang and C. Weng, More on pooling spaces, Discrete
Math. 308 (2008)  6330--6338.

\bibitem{hww} T. Huang, K. Wang and C. Weng, More pooling spaces associated with some
finite geometries, European J.  Combin. 29 (2008) 1483--1491.


\bibitem{HW}
T. Huang and C. Weng, Pooling spaces and non-adaptive pooling designs,
Discrete Math. 282 (2004) 163--169.


\bibitem{Macula}
A. J. Macula, A simple construction of $d$-disjunct matrices with
certain constant weights, Discrete Math. 162 (1996) 311--312.

\bibitem{Macula2}
A. J. Macula, Error-correcting non-adaptive group testing with
 $d^e$-disjunct matrices, Discrete Appl. Math. 80 (1997) 217--222.

\bibitem{NG}
J. Nan and J. Guo, New error-correcting pooling designs associated
with finite vector spaces, J. Comb. Optim. 20 (2010) 96--100.


\bibitem{Ngo}
H. Ngo and D. Du, New constructions of non-adaptive and
error-tolerance pooling designs, Discrete Math. 243 (2002)
167--170.

\end{thebibliography}
\end{document}